\documentclass[letterpaper,twoside,twocolumn]{IEEEtran}
\IEEEoverridecommandlockouts        

\usepackage[utf8]{inputenc}
\usepackage{multicol}
\usepackage{graphicx,float,wrapfig,color}
\usepackage{url}
\usepackage{amssymb,amsmath,amsfonts,layout}
\usepackage{makeidx}
\usepackage{blkarray,multirow}
\usepackage{cite}
\usepackage{enumerate}
\usepackage{algpseudocode}
\usepackage{algorithm}

\usepackage{amsthm}
\usepackage{bbm}
\usepackage{caption}
\usepackage{hyperref}
\usepackage{comment}

\algnewcommand{\Inputs}[1]{%
	\State \textbf{Inputs:}
	\Statex \hspace*{\algorithmicindent}\parbox[t]{.8\linewidth}{\raggedright #1}
}
\algnewcommand{\Initialize}[1]{%
	\State \textbf{Initialize:}
	\Statex \hspace*{\algorithmicindent}\parbox[t]{.8\linewidth}{\raggedright #1}
}
\algnewcommand{\Outputs}[1]{%
	\State \textbf{Outputs:}
	\Statex \hspace*{\algorithmicindent}\parbox[t]{.8\linewidth}{\raggedright #1}
}

\DeclareMathOperator{\A}{\mathcal{A}}
\DeclareMathOperator{\B}{\mathcal{B}}

\DeclareMathOperator{\N}{\mathcal{N}}

\DeclareMathOperator{\R}{\mathbb{R}}

\DeclareMathOperator{\E}{\mathcal{E}}

\newcommand{\argmin}{\operatorname{argmin}}

\newcommand{\setdef}[2]{\{#1 \; | \; #2\}}

\newcommand{\real}{\ensuremath{\mathbb{R}}}

\renewcommand{\bar}{\overline}

\newcommand{\rev}[1]{{\color{blue} #1}}
\renewcommand{\rev}[1]{#1}
\newcommand{\revII}[1]{{\color{blue} #1}}
\renewcommand{\revII}[1]{#1}

\begin{document}

\title{On Hybrid Quantum and Classical Computing Algorithms for Mixed-Integer Programming
\thanks{
The authors are with the National Renewable Energy Laboratory (NREL), \{{chinyao.chang, eric.jones, yiyun.yao, peter.graf, jain.rishabh}\}@nrel.gov. This work was authored in part by the National Renewable Energy Laboratory (NREL), operated by Alliance for Sustainable Energy, LLC, for the U.S. Department of Energy (DOE) under Contract No. DE-AC36-08GO28308. This work was supported by the Laboratory Directed Research and Development (LDRD) Program at NREL. The views expressed in the article do not necessarily represent the views of the DOE or the U.S. Government. The U.S. Government retains and the publisher, by accepting the article for publication, acknowledges that the U.S. Government retains a nonexclusive, paid-up, irrevocable, worldwide license to publish or reproduce the published form of this work, or allow others to do so, for U.S. Government purposes.
}
}
\author{Chin-Yao Chang, Eric Jones, Yiyun Yao, Peter Graf, and Rishabh Jain
}

\maketitle
\begin{abstract}
    Quantum computing is emerging as a new computing resource that could be superior to conventional computing for certain classes of optimization problems. However, in principle, most existing approaches to quantum optimization are intended to solve unconstrained binary programming problems, while mixed-integer linear programming is of most interest in practice. 
    We attempt to bridge the gap between the capability of quantum computing and real-world applications by developing a new approach for mixed-integer programming. The approach applies Benders decomposition to decompose the mixed-integer programming into binary programming and linear programming sub-problems, which are solved by a noisy intermediate-scale quantum processor and conventional processor, respectively. The algorithm is provably able to reach the optimal solution of the original mixed-integer programming problem. The algorithm is tested on a D-Wave 2000Q quantum processing unit and is shown to be effective for small-scaled test cases. \rev{We also test the algorithm on a mixed-integer programming inspired by power system applications. Many insights are drawn from the numerical results for both the capabilities and limitations of the proposed algorithm.}
\end{abstract}

%\begin{keywords}
%Benders decomposition, mixed-integer programming, quantum annealing, quantum computing, optimal power flow.
%\end{keywords}

\section{Introduction}
Recent years have have witnessed exciting progress in the field of quantum computing (QC). Many QC algorithms have been developed to solve various types of problems, optimization being among them. For integer factorization, Shor's algorithm is exponentially faster than the best-known classical algorithms~\cite{shor1994algorithms}. Grover's algorithm \cite{grover1997quantum} provides a quadratic speedup in database search problems. The Harrow-Hassidim-Lloyd (HHL) algorithm~\cite{harrow2009quantum} provides an exponential speedup over classical computing (CC) in solving linear equations. And quantum machine learning~(
\cite{zhang2020recent} and references therein
%\cite{schuld2015introduction,biamonte2017quantum,zhang2020recent} and references therein
), built largely on the aforementioned QC algorithms, is also emerging as a popular research direction. However, many of these quantum algorithms with provable speedups are contingent upon the as-yet unrealized construction of fault-tolerant quantum processors \cite{nielsen2002quantum}. The potential of QC is immense, but the classes of optimization problems that near-term, noisy intermediate-scale quantum (NISQ) processors are able to address remain restricted \cite{preskill2018quantum}. This paper is intended to increase the impact of NISQ optimization by leveraging its capability with classical co-processing to solve mixed-integer programming (MIP) problems, which are of interest in a wide range of real-world applications. 

Many practical control and optimization problems can be formulated as MIPs because the controls (or decision variables in the context of optimization) are a mix of discrete and continuous actions. MIPs occur in various applications such as the classical traveling salesman or routing problems~\cite{matai2010traveling,malandraki1992time}, supply chain planning~\cite{pochet2006production,ozceylan2013mixed}, \revII{energy market operation~\cite{streiffert2005mixed}, and power system expansion planning \cite{alizadeh2011reliability}.}
%and power systems design and operation \cite{streiffert2005mixed,alizadeh2011reliability,zhang2017optimal,li2005price}. 
Finding an optimal solution to a MIP in principle requires complete combinatorial enumeration, leading to exponential complexity. In addition to state-of-the-art branch-and-cut (or branch-and-bound) methods, various specialized algorithms have been developed to efficiently solve some special cases of MIP. This reflects the importance and challenges involved in MIP. 

Quantum processors are able to directly represent exponentially many solutions on linearly scaling hardware resources through qubit state superposition and leverage a uniquely quantum phenomenon known as amplitude amplification in order to soften the hard exponential scaling of combinatorial enumeration. Some variants of the amplitude amplification mechanism are thought to be able to soften the time-complexity scaling exponent for combinatorial optimization so effectively that they outclass the best possible classical heuristics, opening the door to a discussion of quantum advantage.  %\rev{\cite{arute2019quantum,zhong2020quantum}} 
However, most near-term quantum optimization heuristics rely on solving Quadratic Unconstrained Binary Optimization (QUBO) problem formulations. Any optimization problem for which there is an analytic expression can be reformulated as a QUBO, and many formal mathematical and computational problems have been recast in QUBO form, including the minimal dominating set formulation of optimal phasor measurement unit placement on the power grid \cite{jones2020computational}, the Markov decision process \cite{jones2020k}, and many NP difficult problems \cite{lucas2014ising}. With well established CC algorithms for convex optimization, decomposing a MIP problem into a series of one or more QUBO and convex optimization sub-problems, which are respectively solved by QC and CC, is a promising way to leverage QC capabilities. \rev{There are only a few works that engage in this hybrid QC-CC concept.}  In~\cite{gambella2020multi}, an alternating direction method of multipliers (ADMM)-based hybrid QC-CC algorithm was developed to solve a class of MIP with constraints represented as equalities coupling binary and continuous variables. \rev{However, the ADMM-based methods can not guarantee the convergence due to the non-smooth nature of MIP.} In~\cite{ajagekar2020quantum}, hybrid QC-CC algorithms were proposed for job-shop scheduling problems, 
cell formation manufacturing problems, 
and vehicle routing problems. \rev{The latter two problems are non-convex even without considering the discrete binary variables. Therefore, some heuristic methods instead of algorithms with exact ones were developed for them in~\cite{ajagekar2020quantum}.}
% Monte Carlo:\cite{orus2019quantum} \\

In this paper, we propose applying Benders decomposition~\cite{benders1962partitioning} for a hybrid QC-CC algorithm that is provably convergence for a class of MIP. The equivalence between the original MIP and a binary programming problem is derived. The binary programming sub-problem can involve a large number of constraints. In order to avoid this problem, we follow the Benders decomposition iterative procedure by only adding one constraint at a time to gradually tighten the optimization search space. With a sufficient number of back-and-forth iterations of solving the QUBO and linear programming (LP) (used to identify and add constraints to QUBO) sub-problems, the algorithm in principle can find an optimal solution. However, there are many challenges for the practical implementations. For small-scale MIPs, the algorithm can find the optimal solution with 2000Q D-Wave System, while for a more engaging formulation of MIP inspired by power systems applications, the algorithm can only reach a near global optimal solution. Limitations of reformulating inequality constraints to QUBOs and current quantum annealing hardware post some challenges for the practicability of the algorithm, which will be detailed in the numerical studies. 

The rest of the paper is organized as follows. In section~\ref{sec:binary_form}, we begin with a simplified form of MIP in order to review the core ideas of Benders decomposition. Section~\ref{sec:general_form} extends the results in section~\ref{sec:binary_form} for general MIP problems. Section~\ref{sec:reformulation} presents our numerical results.

\textit{Notation:} Let $\R$ and $\R_+$ respectively denote the set of real and non-negative real numbers; $\mathbb{N}$ denotes the set of natural numbers. 
%Given a matrix $A$, $A^\top$ denotes its transpose, and ${A\succ(\prec) 0}$ denotes that $A$ is positive (negative) definite. The matrix $I_n \in \R^{n\times n}$ is the $n \times n$ identity matrix. For $x\in\real^n$, $\norm{x}_2$ denotes its Euclidean norm, $\text{diag}\{x\}$ is a diagonal matrix with the elements of $x$ on the main diagonal. Further, given a set $\mathcal{X}\subset \real^n$, $\text{Proj}_{\mathcal{X}}\{x\}$ denotes the projection of $x$ onto $\mathcal{X}$. Given a $Q\succ 0$, we define $\norm{v}_Q^2 = v^\top Q v$. We use subscript $t$ to denote a vector value at time $t$; and $v_{i,t}$ to denote the $i^{th}$ entry of the vector $v_t$ at time $t$.  
Given a set $\mathcal{S}$, $|\mathcal{S}|$ is the cardinality of it. For a scalar $a\in\real$, $\lfloor a \rfloor$ is the largest integer less than $a$.

\section{Hybrid Algorithm for Simplified Mixed-Integer Programming}\label{sec:binary_form}
We consider a special case of  MIP in this section, shown in the following
\begin{subequations}\label{eq:p1}
\begin{align}
   & \min_{y\in\{0,1\}^{n_y},z\in\real^{n_z}_+}f(y), \label{eq:p1-1} \\
   & \text{s.t. } A_y y + A_z z = b, \label{eq:p1-2} \\
   & \quad\quad Cz \leq d,\label{eq:p1-3}
\end{align}
\end{subequations}
where $n_y$ and $n_z$ are respectively the number of binary and continuous variables, $f:\{0,1\}^{n_y}\mapsto\real$ is a quadratic function (can be convex or non-convex),  $b\in\real^{m_b}$, $d\in\real^{m_d}$, $m_b$ and $m_d$ are the number of equality and inequality constraints, and $A_y$, $A_z$, $C$ are in proper dimensions. We first decompose~\eqref{eq:p1} into two subproblems such that one is a pure binary programming problem that can be solved by QC, and the other can be efficiently solved by conventional computers. Let us rewrite~\eqref{eq:p1} in the following equivalent form
\begin{align}\label{eq:p2}
    & \min_{y\in\{0,1\}^{n_y},}f(y) + q(y),
\end{align}
where
\begin{align}\label{eq:qfunction}
    q(y) = \min_{z\in\real^{n_z}_+}{0}, \text{ s.t. } \eqref{eq:p1-2}-\eqref{eq:p1-3} \text{ hold}.
\end{align}
We can view~\eqref{eq:qfunction} as a feasibility problem. In the context of optimization theory~\cite{boyd2004convex},  primal optimization~\eqref{eq:qfunction} has the associated dual problem formulation. Let $\alpha\in\real^{m_b}$, $\beta\in\real_+^{m_d}$, $ \lambda\in\real_+^{n_z}$ be the dual variables associated with constraints~\eqref{eq:p1-2}, \eqref{eq:p1-3} and $z\in\real^{n_z}_+$ respectively. The dual problem of~\eqref{eq:qfunction} is given as
\begin{align*}
    \max_{\alpha\in\real^{m_b}, \beta\in\real_+^{m_d}, \lambda\in\real_+^{n_z}} & \bigg( \inf_{z} \Big(\alpha^\top (A_y y + A_z z - b) \\
    &\hspace{7mm} + \beta^\top (Cz -d) - \lambda^\top z \Big)\bigg),
\end{align*}
which is equivalent to the following
\begin{align}
    \max_{\alpha\in\real^{m_b}, \beta\in\real_+^{m_d}, \lambda\in\real_+^{n_z}} & \bigg( \inf_{z} \Big((\alpha^\top A_z+ \beta^\top C - \lambda^\top)z\Big) \label{eq:qfunction-d} \\ \nonumber
    &\hspace{7mm}  +  \alpha^\top(A_y y - b) - \beta^\top d \Bigg). 
\end{align}
\rev{Because Slater's condition holds for the linear programming~\eqref{eq:qfunction}, 
%\footnote{Slater's condition holds for most linear programming problems.}
the duality gap between the primal problem,~\eqref{eq:qfunction}, and dual problem,~\eqref{eq:qfunction-d}, is zero.} This implies that the optimal value of~\eqref{eq:qfunction-d} should be zero if~\eqref{eq:qfunction} is feasible. On the other hand, if~\eqref{eq:qfunction} is infeasible, then the optimal value of~\eqref{eq:qfunction-d} is infinity. The optimal value of~\eqref{eq:qfunction-d} should be non-negative either way. Optimization~\eqref{eq:qfunction-d} has the optimal value non-negative only if $\alpha^\top A_z + \beta^\top C - \lambda^\top = 0$. Otherwise, there exists a $z$ such that~\eqref{eq:qfunction-d} has the optimal value $-\infty$. We can then rewrite~\eqref{eq:qfunction-d} in the following equivalent form
\begin{align}\label{eq:qfunction-d2}
    \max_{\alpha\in\real^{m_b}, \beta\in\real_+^{m_d}}& \alpha^\top(A_y y - b) - \beta^\top d,  \\ \nonumber
    \text{ s.t. }& \alpha^\top A_z + \beta^\top C \geq 0.
\end{align}
The optimal value of~\eqref{eq:qfunction-d2} is either zero or infinity, which respectively indicates~\eqref{eq:qfunction} is feasible or infeasible for the given $y$. The $y$ such that~\eqref{eq:qfunction} is feasible should make the optimum of~\eqref{eq:qfunction-d2} zero. Optimization~\eqref{eq:qfunction-d2} with $y$ given has the optimum being zero if and only if 
\begin{align}\label{eq:extreme_rays}
    \alpha^\top_k(A_y y - b) - \beta^\top_k d \leq 0,\quad \forall (\alpha_k, \beta_k) \in \Gamma,
\end{align}
where $\Gamma$ is the set of extreme rays of the set $\setdef{(\alpha, \beta)}{\alpha^\top A_z + \beta^\top C \geq 0, \; \beta \geq 0}$. To this end, we can conclude that \eqref{eq:qfunction} is feasible if and only if \eqref{eq:extreme_rays} holds. Replacing $q(y)$ in~\eqref{eq:p2} by~\eqref{eq:extreme_rays} gives an equivalent formulation of~\eqref{eq:p1}, shown as follow
\begin{align}\label{eq:p3}
   & \min_{y\in\{0,1\}^{n_y}}f(y), \text{ s.t. \eqref{eq:extreme_rays} holds}.
\end{align}
Optimization~\eqref{eq:p3} is a \rev{quadratic} binary programming problem and in principle can be solved by QC. The potential issue is that the number of extreme rays, or $|\Gamma|$, may incur too many inequality constraints, which in turn would require a prohibitive number of qubits. 

In the context of Benders decomposition, optimization~\eqref{eq:p3} without constraint~\eqref{eq:extreme_rays} is known as the master problem, and~\eqref{eq:extreme_rays} collects the set of ``cuts'' (or constraints) that may gradually add to the master problem. If the optimal solution $y_0$ of the unconstrained~\eqref{eq:p3} has~\eqref{eq:qfunction-d2} unbounded, it only indicates that $y_0$ is not a feasible solution, but it does not provide information on which cut to add. For convenience of determining which cut to add to the master problem, we define~\eqref{eq:qfunction-d3}, which is the same as~\eqref{eq:qfunction-d2} except that bounds on $\alpha$ and $\beta$ are added.
\begin{align}\label{eq:qfunction-d3}
    \max_{|\alpha| \leq \bar{\alpha}, \; |\beta| \leq \bar{\beta}, \beta\geq 0}& \alpha^\top(A_y y - b) - \beta^\top d, \\ \nonumber 
    \text{ s.t. } & \alpha^\top A_z + \beta^\top C \geq 0, 
\end{align}
where $\bar{\alpha}\in\real^{m_b}_+$ and $\bar{\beta}\in\real^{m_d}_+$. Both $\bar{\alpha}$ and $\bar{\beta}$ are finite. The optimal solution of~\eqref{eq:qfunction-d3} is clearly bounded. In addition, the optimal solution of~\eqref{eq:qfunction-d3} should be one of the extreme rays in $\Gamma$. Repeatedly solving~\eqref{eq:qfunction-d3} leads to a sequence of cuts to add to unconstrained~\eqref{eq:p3}. Algorithm~\ref{algo:decomposition} summarizes the steps described above.
\begin{algorithm}[H]
\caption{}
\begin{algorithmic}[1] 
\State \textbf{Initialize} $\Xi =\emptyset$, $t= 0$, $\xi = 0$
\While {$\xi =0$}
    \State Find $y_t = \argmin_{y\in\{0,1\}^{n_y}}f(y), \text{ s.t. } y\in\Xi$
    \State Find $(\alpha_t,\beta_t)$ by solving~\eqref{eq:qfunction-d3} for the given $y_t$
    \If {$\alpha^\top_t(A_y y_t - b) - \beta^\top_t d > 0$} 
        \State $\Xi = \Xi \cup \setdef{y}{\alpha^\top_t(A_y y - b) - \beta^\top_t d  \leq 0}$
        \State $t\mapsto t+1$
    \Else 
        \State Conclude $y_t$ is the optimal solution
        \State Solve~\eqref{eq:qfunction} with given $y_t$ Find the optimal $z$
        \State {Set $\xi = 1$} 
    \EndIf
\EndWhile
\end{algorithmic}
\label{algo:decomposition}
\end{algorithm}
Note that in the worst case, one would need to run $|\Gamma|$ iterations before finding the optimal solution. One alternative may be setting an upper bound on the number of iterations. If Algorithm~\ref{algo:decomposition} is unable to find the optimal solution within the maximal number of iterations, one can assume that the incumbent solution is close to the global optimum and use that to assist the conventional computer to solve~\eqref{eq:p1}.
% \begin{align}\label{eq:qfunction}
 %   & q(y) = \begin{cases} \inf & \text{if }\not\exists z\in\real^{n_z}_+ \text{ s.t. }~\eqref{eq:p1-2}-\eqref{eq:p1-3} \text{ hold},\\
 %   0 & \text{if }\exists z\in\real^{n_z}_+ \text{ s.t. }~\eqref{eq:p1-2}-\eqref{eq:p1-3} \text{ hold}.
 %   \end{cases}
%\end{align}

\section{General Mixed-Integer Programming}
\label{sec:general_form}
We now generalize~\eqref{eq:p1} by adding terms linear in the continuous variables to the objective function. The equalities in~\eqref{eq:p1-2} are replaced by inequalities.
\begin{subequations}\label{eq:p1g}
\begin{align}
   & \min_{y\in\{0,1\}^{n_y},z\in\real^{n_z}_+}f(y) + g^{\top}z, \label{eq:p1g-1} \\
   & \text{s.t. } A_y y + A_z z \leq  b, \label{eq:p1g-2} \\
   & \quad\quad Cz \leq d,\label{eq:p1g-3}
\end{align}
\end{subequations}
Following the steps in the last section, we can rewrite optimization~\eqref{eq:p1g} in the following
\begin{comment}
\begin{subequations} \label{eq:equivalent_dual}
\begin{align} 
    &\min_{q\in\real} q, \label{eq:equivalent_dual_cost} \\
    & \text{ s.t. } \alpha^\top_k(A_y y - b) - \beta^\top_k d \leq 0,\quad \forall (\alpha_k, \beta_k) \in \Gamma_r, \label{eq:equivalent_dual_ray} \\ 
    & \quad\quad \alpha^\top_k(A_y y - b) - \beta^\top_k d \leq q,\quad \forall (\alpha_k, \beta_k) \in \Gamma_p. \label{eq:equivalent_dual_point}
\end{align}
\end{subequations}
Combining~\eqref{eq:p2g} and~\eqref{eq:equivalent_dual} gives
\end{comment} 
\begin{subequations} \label{eq:pg3}
\begin{align}
   & \min_{y\in\{0,1\}^{n_y},q\in\real}f(y) + q, \label{eq:pg3-1} \\
   %\text{ s.t. \eqref{eq:equivalent_dual_ray} and \eqref{eq:equivalent_dual_point} hold}.
   & \text{ s.t. } \alpha^\top_k(A_y y - b) - \beta^\top_k d \leq 0,\quad \forall (\alpha_k, \beta_k) \in \Gamma_r, \label{eq:pg3_ray} \\ 
    & \quad\quad \alpha^\top_k(A_y y - b) - \beta^\top_k d \leq q,\quad \forall (\alpha_k, \beta_k) \in \Gamma_p, \label{eq:pg3_point}
\end{align}
\end{subequations}
where $\Gamma_r$ and $\Gamma_p$ respectively are the set of extreme rays and extreme points of the set $\setdef{(\alpha, \beta)}{\alpha^\top A_z + \beta^\top C +g^{\top} \geq 0, \; \alpha \geq 0, \; \beta \geq 0}$. Different from~\eqref{eq:p3}, \eqref{eq:pg3} is not a binary programming because it has one continuous scalar variable $q$. This is a trade-off that we need to make when we consider a more general class of MIP.

We follow a similar Benders decomposition process to solve~\eqref{eq:pg3}. Consider the following dual optimization (similar logic to \eqref{eq:qfunction-d3})
\begin{align}\label{eq:qfunction_g-d2_bdd}
    \max_{|\alpha| \leq \bar{\alpha}, \; |\beta| \leq \bar{\beta}, \beta\geq 0} & \alpha^\top(A_y y - b) - \beta^\top d,  \\ \nonumber 
    \text{ s.t. }& \alpha^\top A_z + \beta^\top C + g^{\top} \geq 0.
\end{align}
The optimal solution of~\eqref{eq:qfunction_g-d2_bdd} should be in the set $\Gamma_r$ or $\Gamma_p$. Similar to the previous section, repeatedly solving~\eqref{eq:qfunction_g-d2_bdd} leads to a sequence of cuts to add to the unconstrained~\eqref{eq:pg3}. However, any added cut depends on whether the solution of~\eqref{eq:qfunction_g-d2_bdd}, $(\alpha^{\star},\beta^{\star})$, is in $\Gamma_r$ or $\Gamma_p$. We use the following simple logic to determine which set the solution is in
%\begin{subequations} 
\begin{align} \label{eq:feas_opt}
%& \text{If $\exists k$ s.t. $|\alpha^\star(k)| =  \bar{\alpha}^\star(k)$, or $|\beta^\star(k)| = \bar{\beta}^\star(k)$}, \\
%&\quad (\alpha^{\star},\beta^{\star}) \in \Gamma_r, \\
%&\text{otherwise,} \\
%&\quad  (\alpha^{\star},\beta^{\star}) \in \Gamma_p. 
\begin{cases}
    (\alpha^{\star},\beta^{\star}) \in \Gamma_r & \text{if $\exists k$ s.t. $|\alpha^\star(k)| = \bar{\alpha}^\star(k)$}, \\
    & \quad\quad\quad \text{or $|\beta^\star(k)| = \bar{\beta}^\star(k)$}, \\
    (\alpha^{\star},\beta^{\star}) \in \Gamma_p & \text{otherwise},
\end{cases}
\end{align}
%\end{subequations}
where $\alpha^{\star}(k)$ is the $k^{th}$ entry of $\alpha^{\star}$ and similar notation applies for other vector variables.
If an entry of $\alpha^{\star}$ (or $\beta^{\star}$) reaches the upper bound, then it is in $\Gamma_r$ if $\bar{\alpha}$ and $\bar{\beta}$ are chosen large enough, because the optimal points in $\Gamma_r$ are unbounded if we do not add the bounds $\bar{\alpha}$ and $\bar{\beta}$. On the contrary, the extreme points are in the bounded set and constitute a strict subset of $\setdef{(\alpha,\beta)}{|\alpha| < \bar{\alpha}, \; |\beta| < \bar{\beta} }$. The logic shown in~\eqref{eq:feas_opt} is accurate as long as we have $\bar{\alpha}$ and $\bar{\beta}$ large enough. 

Algorithm~\ref{algo:decomposition-2} summarizes the steps for finding the optimal point of~\eqref{eq:p1g}. The main difference between Algorithm~\ref{algo:decomposition} and~\ref{algo:decomposition-2} lies in the form of the cuts. The cuts added in Algorithm~\ref{algo:decomposition} are known as ``feasibility cuts'' because those cuts are added for the purpose of enforcing a feasible solution for the original optimization~\eqref{eq:p1}.  The intuition behind the feasibility cut is that the continuous variables, $z$, only affect the feasibility of the optimization because they only appear in the constraint of~\eqref{eq:p1}.  Optimization~\eqref{eq:p1g} that  Algorithm~\ref{algo:decomposition-2} targets to solve has $z$ simultaneously in the constraint set and objective function. As a result, in the course of iterative steps of Algorithm~\ref{algo:decomposition-2}, we may add feasibility cuts (similar to Algorithm~\ref{algo:decomposition}) or ``optimality cuts'' to reflect  appearance of $z$ in the constraint set and the objective function.  Feasibility and optimality cuts are very similar, the only difference between them being the dependency on the variable $q$ as shown in line 10 and 13 of Algorithm~\ref{algo:decomposition-2}.

%Both Algorithms~\ref{algo:decomposition} and~\ref{algo:decomposition-2} up to this point can only be applied to mixed-integer linear programming (MILP). However, we can generalize the results to the following mixed-integer quadratic programming by using linear matrix inequalities (LMIs) techniques
%\begin{subequations}\label{eq:p1g2}
%\begin{align}
%   & \min_{y\in\{0,1\}^{n_y},z\in\real^{n_z}_+}f(y) + z^{\top}Gz g^{\top}z, \label{eq:p1g2-1} \\
   %& \text{s.t. } A_y y + A_z z \leq  b, \label{eq:p1g2-2} \\
   %& \quad\quad Cz \leq d,\label{eq:p1g2-3}
%\end{align}
%\end{subequations}

\begin{algorithm}[]
\caption{}
\begin{algorithmic}[1] 
\State \textbf{Initialize} $\Xi =\emptyset$, $t= 0$, $\xi = 0, \rev{\gamma} = 0, s_{pt} = 0, s_{dt} = 0$
\While {$\xi =0$}
    \State Solve the optimization
    \begin{align*}
    (y_t,p_t) = &\argmin_{y\in\{0,1\}^{n_y},q\in\real}f(y) + \rev{\gamma} q, \\ \quad\quad\text{ s.t. } &(y,q)\in\Xi
    \end{align*}
    \State Find $(\alpha_t,\beta_t)$ by solving~\eqref{eq:qfunction_g-d2_bdd} for the given $y_t$
    \State Extract the optimal value of the primal problem: 
    \begin{align*}
    s_{pt} = f(y_t) + \rev{\gamma} q_t 
    \end{align*}
    \State Extract the optimal value of~\eqref{eq:qfunction_g-d2_bdd} given as 
    \begin{align*}
    s_{dt} = \alpha^\top_t(A_y y_t - b) - \beta^\top_t
    \end{align*}
    \State Use~\eqref{eq:feas_opt} to identify whether $(\alpha_t,\beta_t)$ is in $\Gamma_r$ or $\Gamma_p$
    \If {$s_{pt} < s_{dt}$} 
        \If {$(\alpha_t,\beta_t)\in\Gamma_r$}
            \State $\Xi = \Xi \cup \setdef{(y,q)}{\alpha^\top_t(A_y y - b) - \beta^\top_t d  \leq 0}$
            \State $t\mapsto t+1$    
        \Else
            \State $\Xi = \Xi \cup \setdef{(y,q)}{\alpha^\top_t(A_y y - b) - \beta^\top_t d  \leq q}$
            \State Set $\rev{\gamma} = 1$
            \State $t\mapsto t+1$    
        \EndIf
    \Else 
        \State Conclude $(y_t,q_t)$ is the optimal solution
        \State Find the optimal $z$ with the given $y_t$ 
        \State {Set $\xi = 1$} 
    \EndIf
\EndWhile
\end{algorithmic}
\label{algo:decomposition-2}
\end{algorithm}

\section{Reformulation of Binary Programming} \label{sec:reformulation}
Both Algorithms~\ref{algo:decomposition} and~\ref{algo:decomposition-2} pose a binary programming problem to be solved via quantum computing. The binary programming problem becomes constrained  when cuts are added during each of the while loop iterations. However, because most existing quantum optimization heuristics are designed only to accept QUBO problems, we need to recast hard constraints as soft, penalty functions added to the objective function. We will focus on a penalization method for Algorithm~\ref{algo:decomposition} for simplicity. The logic of the penalization method for Algorithm~\ref{algo:decomposition-2} is very similar.

In the following we write down  in explicit form the optimization in line 3 of Algorithm~\ref{algo:decomposition}
\begin{align}\label{eq:binary_qc0}
    \min_{y\in \{0,1\}^{n_y} } f(y), \quad \text{ s.t. } \hat{A}_t y - \hat{b}_t \leq 0, 
\end{align}
where $\hat{A}_t \in \real^{t\times n_y}$ and $\hat{b}_t\in \real^{t}$. The number of rows for both $\hat{A}_t$ and $\hat{b}_t$ increases by one whenever $t$ increases, which can be equivalently interpreted as an added cut. We consider penalizing each cut by a quadratic objective function shown in the following
\begin{subequations}\label{eq:binary_qc1}
\begin{align}
    &\min_{y\in \{0,1\}^{n_y}, x_k\in \{0,1\}^{n_{x_k}}, k = 1,\cdots,t} f(y) \label{eq:binary_qc1-1} \\ 
    &+ w_t\sum_{k = 1}^{t} \Big(a_{k} y - \hat{b}_t(k) + w_{x_k} \sum_{j = 1}^{n_{x_k}}2^{(j-1)}x_k(j)\Big)^2, \label{eq:binary_qc1-2}
\end{align}
\end{subequations}
where $a_k \in \real^{n_y}$ is the $k^{th}$ row of the matrix $\hat{A}_t$, $w_t\in\real_+$ and $w_{x_k}\in\real_+$ are penalty weights, and $x_k(j)$ is the $j^{th}$ element of $x_k$. The ideal penalty weights $w_t$ and $w_{x_k}$ have the following properties: (i) for any $y$ such that the inequalities in~\eqref{eq:binary_qc0} hold,  there exists feasible $x_k$, $k=1,\cdots,t$, such that \eqref{eq:binary_qc1-2} is zero; (ii) all solutions violating the constraint~\eqref{eq:binary_qc0} are not the optimal $y$ for~\eqref{eq:binary_qc1}. If (i) and (ii) hold true, then \eqref{eq:binary_qc1} and \eqref{eq:binary_qc0} are equivalent in the sense that they have the same optimal $y$.

We first analyze condition (i). If $a_k$ and $\hat{b}_t(k)$ have integer values, then we can always choose $w_{x_k}=1$ for all $k=1,\cdots,t$ to ensure that for any $y$ feasible for~\eqref{eq:binary_qc0}, there exists a $x_k$ such that $a_k y - \hat{b}_t(k) + w_{x_k} \sum_{j = 1}^{n_{x_k}}2^{(j-1)}x_k(j) = 0$ as long as $n_{x_k}$ is sufficiently large. When the values of $a_k$ and $\hat{b}_t(k)$ are not integer, we choose a small value for $w_{x_k}$ so that there exists a $x_k$ such that $a_k y - \hat{b}_t(k) + w_{x_k} \sum_{j = 1}^{n_{x_k}}2^{(j-1)}x_k(j) = 0$ for any feasible $y$. In this case, we have a trade-off between the accuracy and the number of qubits. If $w_{x_k}$ is chosen very small for accuracy purposes, a large $n_{x_k}$ is necessary, which implies a large number of qubits. Currently we choose $w_{x_k}$ in heuristic ways, but a rigorous way of choosing $w_{x_k}$ may be attainable. For condition (ii), there always exists a large enough $w_t$ so that (ii) holds. However, in the context of QC, if $w_t$ is overwhelmingly large compared to the actual objective function $f(y)$, then noise could impair the QC from reaching the ground state, i.e., optimal solution. As a result, QC may reach a solution that satisfies the penalty constraints while $f(y)$ is not minimized. Again we compromise by choosing $w_t$ in heuristic ways.

%The term $w_{x_k} \sum_{j = 1}^{n_{x_k}}2^{(j-1)}x_k(j)$ is always non-negative, which implies to minimize the cost of~\eqref{eq:binary_qc1} other than $f(y)$, $y$ should be chosen such that $a_k y - \hat{b}_t(k) \leq 0$ for all $k=1,\cdots,t$. 

\section{Numerical Studies}\label{sec:simulation}
\rev{There are many potential applications that can be formulated as~\eqref{eq:p1} or~\eqref{eq:p1g} so that the proposed algorithm can become useful. For example, planning of transmission systems, co-optimization of grid and transportation. In this section, we consider two different problem formulations for the algorithm validation purpose.} One is a small-scale MIP and the other is a more engaging and realistic MIP for a unit commitment problem for power systems. For the small-scale example, we implement the proposed algorithm on an actual D-Wave system while for the unit commitment problem, we use a simulated annealing sampler as an alternative to the D-Wave system because of current size and connectivity limitations.

\subsection{Small-Scale Test Case for D-Wave System}\label{sec:simulation-1}
We validate the proposed algorithm by running small-scale MIP optimizations on a 2000Q D-Wave quantum processing unit (QPU). The 2000Q D-Wave QPU is constructed from a 2,048 qubit $(16,4)$-chimera hardware graph topology consisting of a $16\times 16$ array of $8-$qubit chimera unit cells, each of which is organized as a $K_{4,4}$ complete bipartite graph. More details can be found in~\cite{junger2019performance}. We accessed the D-Wave system by their ${Leap}^2$ project.

The objective function for the binary variables takes the linear form $f(y) = g_y^\top y$. 
\begin{subequations}\label{eq:MILP_num}
\begin{align}
    & g_y = \begin{bmatrix} 8.75 & 4.95 & 6.2 & 8.4 \end{bmatrix}^{\top}, \\ 
    & A_y = \begin{bmatrix} 5 & 3 & 4 & 6 \\ 2.5 & 1.2 & 2 & 1.8 \\ 1.5 & 0.9 & 1.6 & 2.4
    \end{bmatrix}, \\
    & A_z = \begin{bmatrix} 1 & 1 & 1 & 1 \\ 0.8 & 0.7 & 0.6 & 0.3 \\ 0.6 & 0.7 & 0.8 & 0.9
    \end{bmatrix}, \\
    & b = \begin{bmatrix} 25 & 12.5 & 12.5 \end{bmatrix}^{\top}.
\end{align}
\end{subequations}
Both $C$ and $d$ are zero in this test case. Using our proposed decomposition approach, the original MIP is decomposed into a binary programming and LP, which are solved in the alternating sequence described in Algorithm~\ref{algo:decomposition}. The D-Wave system solves binary programming and classical computing solves LP. We set the penalty strength ($w_{t}$) at $50$, granularity constant ($w_{x_k}$) at $0.01$, and $n_{x_k} = 10$ for every $t$ and $k$. In addition to those parameters, one needs to define the chain strength for the D-Wave system, which gauges the degree to which multiple physical qubits act as a single logical qubit. On the D-Wave system, the optimization problem is embedded onto its physical Chimera hardware graph. The embedding process may group multiple qubits together to represent one logical variable in the QUBO. Chain strength controls the coupling strength between the grouped qubits. In our simulation setup, we set the chain strength at $75$ for all $t$ to ensure it outweighs the penalty strength. Under this numerical setup,  Algorithm~\ref{algo:decomposition} finds the optimal solution in two iterations. 
%Note that in this case we need to set the penalty strength at a relatively large value, which is mainly due to the small values involved in the second and third rows of $A_y$, $A_z$, and $b$. Those small values force us to have a large penalty strength to outweigh the constraint satisfaction over objective minimization $g_y^\top y$. Ideally the rows that define the linear constraints are at a similar scale. To see this, in our second trial, we scale up the second and third rows of $A_y$, $A_z$, and $b$ by a factor of 10. In which case, we can get the same results with a much smaller penalty strength and chain strength at $50$ and $75$, respectively.

We next simulate a more general MIP with an additional term on the objective function shown in~\eqref{eq:p1g} with $g = [20,2,1,4]^\top$. The remaining constraints are the same as~\eqref{eq:MILP_num}. We use 30 qubits to approximate the continuous scalar $q$ for QC (Line 3 in Algorithm~\ref{algo:decomposition-2}), shown in the following
\begin{align}\label{eq:approx_q}
    q = w_q\sum_{j=1}^{30} 2^{(j-1)}x(j), 
\end{align}
where $w_q=0.01$ and $x(j)\in\{0,1\}$ for all $j=1,\cdots,30$. We set penalty strength and chain strength respectively at $2$ and $3$ in this case. Algorithm~\ref{algo:decomposition-2} reaches the optimal solution in four iterations. Note that the penalty and chain strengths are much smaller in our second test case compared to the first one. The reason is subtle and is related to the upper bounds on $\alpha$ in~\eqref{eq:qfunction-d3} and~\eqref{eq:feas_opt}\footnote{No $\beta$ is defined for the numerical examples because we do not impose any linear inequality constraint~\eqref{eq:p1-3} or~\eqref{eq:p1g-3} in the numerical study.}. For Algorithm~\ref{algo:decomposition}, we set every entry of $\bar{\alpha}$ at 1 and can always decide whether to add a cut or not by checking whether the objective value of~\eqref{eq:qfunction-d3} is strictly positive. However, for Algorithm~\ref{algo:decomposition-2}, we need to differentiate feasibility and optimality cuts, in this case, entries of $\bar{\alpha}$ are set at some arbitrary large values. The large $\bar{\alpha}$ are translated to large coefficients for the feasibility cut defined in Line 10 of Algorithm~\ref{algo:decomposition-2}. To prevent the penalty functions for the feasibility cuts from outweighing  the objective function, we normalize all the cuts in Algorithm~\ref{algo:decomposition-2} by dividing  their coefficients by the square root of the optimal value of~\eqref{eq:feas_opt}. Even with the presence of the normalization, the coefficients of these cuts are still relatively large compared to our first test case, so we choose much smaller penalty and chain strengths. 
%Systematic tuning of the penalty and chain strengths is among our future works.

\subsection{Unit Commitment Test Case}
In this section, we test Algorithm~\ref{algo:decomposition} on a unit commitment problem for power systems applications. The formulation is similar to standard optimal power flow (OPF) with the main difference being in the additional option for every generator $i \in\N_G$ to switch ``on/off'' based on the total power demand, where $\N_G$ is the set of generators. Those switching decisions naturally introduce binary variables. The other main difference from the standard OPF is that we consider planning for a certain time horizon, $T$ hours, with a sampling time of one hour ($\tau = 1,\cdots,T$). The switching between consecutive timestamps is constrained by the minimal locked time of the generators. That is, every generator has to stay in ``on'' (or ``off'') state before it is allowed to change state. Without loss of generality, we assume that the minimal ``on'' (``off'') time, $T_{on}\in\mathbb{N}$ ($T_{off}\in\mathbb{N}$), is strictly larger than one. Otherwise, the minimal ``on'' (``off'') time is not binding for the unit commitment problem. Capturing the locked time constraints introduces additional binary variables. Such a unit commitment problem can be either formulated in the form of~\eqref{eq:p1} or~\eqref{eq:p1g}.

\begin{figure}
    \centering
    \includegraphics{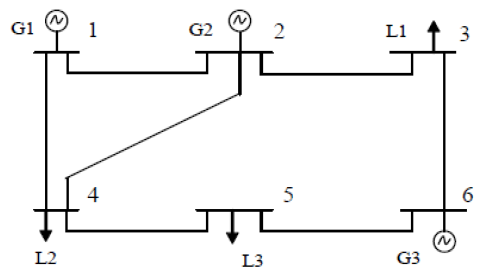}
    \caption{The 6 buses power network with three generators.}
    \label{fig:6bus_network}
\end{figure}

\subsubsection{Formulation for a 6-Bus System}
Given a 6-bus system shown in Figure~\ref{fig:6bus_network}, our goal is to minimize the generation cost subject to the OPF constraints. In the following we go through all the OPF constraints and reformulate them to the ones for the unit commitment problem (with additional binary decisions). First of all, every generator has the upper and lower bounds on the generation limits given as
\begin{align}\label{eq:gen_bdds_ori}
    & \underline{P}_{i} \leq p_i(\tau) \leq \bar{P}_{i}, \quad \forall i \in\N_G, 
\end{align}
for all $\tau = 1,\cdots,T$, where $p_i$ is the active power generation of generator $i$, $\underline{P}_{i}\in\real_+$ and $\bar{P}_i\in\real_+$ respectively are the lower and upper bounds of the generation. Define $y_{G,i}(\cdot) \in \{0,1\}$ for every $i\in\N_G$ such that $y_{G,i}(\cdot)=0$ corresponds to generator $i$ being in the ``on'' state ($y_{G,i}(\cdot)=1$ for the ``off'' state). The constraint above is reformulated as 
\begin{subequations}
\label{eq:gen_bdds}
\begin{align}
    & p_i(\tau) -M\cdot y_{G,i}(\tau) \leq  \bar{P}_{i},  \quad \forall i \in\N_G, \label{eq:gen_bdds-1} \\
    & -\underline{P}_{i} \leq  p_i(\tau) +M\cdot y_{G,i}(\tau),  \quad\forall i \in\N_G,  \label{eq:gen_bdds-2} \\
    & p_i(\tau) + \bar{P}_{i} \cdot y_{G,i}(\tau) \leq  \bar{P}_{i}, \quad \forall i \in\N_G, \label{eq:gen_bdds-3}
\end{align}
\end{subequations}
for all $\tau = 1,\cdots,T$, where $M$ is an arbitrary scalar larger than $\max\{\bar{P}_1,\cdots,\bar{P}_{|\N_G|}\}$. Equation~\eqref{eq:gen_bdds-1} and~\eqref{eq:gen_bdds-2} applies well-known Big M method so that~\eqref{eq:gen_bdds_ori} is active only when $y_{G,i}(\cdot) = 0$. Equation~\eqref{eq:gen_bdds-3} ensures that $p_i(\tau) = 0$ if $y_{G,i}(\tau) = 1$. We next apply similar methods to reformulate the ramp constraints 
\begin{align}
\label{eq:ramp_const_ori}
|p_i(\tau-1) - p_i(\tau)|\leq \bar{P}_{r,i},
\end{align}
where $\bar{P}_{r,i}$ is the maximal ramp rate. The reformulated~\eqref{eq:ramp_const_ori} is given as
\begin{subequations}
\label{eq:ramp_const}
\begin{align}
    & p_i(\tau) - p_i(\tau-1) - M_r \cdot y_{on,i}(\tau-1)\leq   \bar{P}_{r,i}, \label{eq:ramp_const-1} \\
    & p_i(\tau-1) - p_i(\tau) - M_r \cdot y_{off,i}(\tau-1)\leq  \bar{P}_{r,i}, \label{eq:ramp_const-2}
    \end{align}
\end{subequations}
for all $\tau = 2,\cdots,T$, where $y_{on,i}(\cdot)\in\{0,1\}$ with $y_{on,i}(\tau-1) = 1$ indicating that generator $i$ should stay in the ``on'' state at time $\tau$, $y_{off,i}(\cdot)\in\{0,1\}$ with $y_{off,i}(\tau-1) = 1$ indicating that generator $i$ should stay in the ``off'' state at time $\tau$. Equation~\eqref{eq:ramp_const} has the ramp constraint~\eqref{eq:ramp_const_ori} active if $y_{on,i}(\tau-1) = y_{off,i}(\tau-1) = 0$. We will next post constraints such that $y_{G,i}(\tau) = y_{G,i}(\tau-1) = 0$ if $y_{on,i}(\tau-1) = 1$ and $y_{G,i}(\tau) = y_{G,i}(\tau-1) = 1$ if $y_{off,i}(\tau-1) = 1$, so no need to have the ramp constraints when $y_{on,i}(\tau-1) = 1$ or $y_{off,i}(\tau-1) = 1$. For every $i\in\N_G$, $y_{on,i}(\cdot)$, $y_{off,i}(\cdot)$, $y_{G,i}(\cdot)$ are coupled by the following constraints:
\begin{subequations}
\label{eq:OnOff_ind}
\begin{align} 
    & y_{on,i}(\tau-1) + y_{G,i}(\tau) \leq 1, \label{eq:OnOff_ind-1} \\
    & y_{G,i}(\tau) -y_{off,i}(\tau-1) \geq 0, \label{eq:OnOff_ind-2}\\
    & 0 \leq y_{G,i}(\tau)- y_{G,i}(\tau-1) +  y_{on,i}(\tau) \leq 1, \label{eq:OnOff_ind-3} \\
    & 0 \leq y_{G,i}(\tau-1)- y_{G,i}(\tau) +  y_{off,i}(\tau) \leq 1, \label{eq:OnOff_ind-4}  \\
    & \sum_{k = \tau}^{\hat{T}_{on}(\tau) + \tau} y_{on}(k) \geq \Big(y_{on}(\tau) - y_{on}(\tau-1)\Big)\hat{T}_{on}(\tau), \label{eq:OnOff_ind-5}\\
    & \hspace{-3mm}\sum_{k = \tau}^{\hat{T}_{off}(\tau) + \tau} \hspace{-3mm}y_{off}(k) \geq \Big(y_{off}(\tau) - y_{off}(\tau-1)\Big)\hat{T}_{off}(\tau), \label{eq:OnOff_ind-6}
\end{align}
\end{subequations}
where $\hat{T}_{on}(\tau) = T_{on} - \max\{0,T_{on}+\tau - T\}$ and $\hat{T}_{off}(\tau) = T_{off} - \max\{0,T_{off}+\tau - T\}$. Equation~\eqref{eq:OnOff_ind-1} keeps $y_{G,i}(\tau)$ at the ``on'' state if $y_{on,i}(\tau-1) = 1$. Similar logic applies to \eqref{eq:OnOff_ind-2}. Equations~\eqref{eq:OnOff_ind-3} and~\eqref{eq:OnOff_ind-4} respectively have $y_{on,i}(\tau)=1$ and $y_{off,i}(\tau)=1$ when generator $i$ switches at time $\tau$. Equations~\eqref{eq:OnOff_ind-5} and~\eqref{eq:OnOff_ind-6} respectively keep $y_{on,i}(\tau)$ and $y_{off,i}(\tau)$ stay at $1$ before the associated generator is allowed to switch.

All the constraints above are about nodal power injections. The correlations between the buses (nodes) are defined by the graph of the network and the power flow equations. Define $\N$ and $\E$ respectively as the set of buses (nodes) and transmission lines (edges). Constraints associated with the line power flows are given as
\begin{subequations}
\label{eq:powerflow_UC0}
\begin{align}
    & |p_{ik}(\tau)|\leq \bar{P}_{ik}, \forall \{i,k\}\in\E, \label{eq:powerflow_UC0-1} \\
    & P(\tau) = \A P_L(\tau), \label{eq:powerflow_UC0-2} \\
    & P_L(\tau) = \B\A^\top\Theta(\tau),  \label{eq:powerflow_UC0-3}
\end{align}
\end{subequations}
for all $\tau = 1,\cdots,T$, where $p_{ik}(\cdot)$ is the active line power flow for line $\{i,k\}\in\E$, $\bar{P}_{ik}(\cdot)$ is the upper bound of it, $\A$ is the incidence matrix associated with the power network graph (with an arbitrary direction for every line), $\B$ is a diagonal matrix with each of the diagonal element being associated with the impedance of a line, $P(\cdot)\in\real^{|\N|}$ collects the nodal active power injections and $P_L(\cdot)\in\real^{|\E|}$ collects $p_{ik}$ for all $\{i,k\}\in\E$, $\Theta(\cdot)\in\real^{|\N|}$ collects the phase angle of all the buses. Equation~\eqref{eq:powerflow_UC0-3} corresponds to DC power flow formulation which assumes purely inductive lines so the power flow equations can be approximated in the linear form. \rev{More details about the DC power flow can be found in~\cite{van2014dc}.} Because our formulation requires all the continuous variables to be non-negative, we change all the constraints in~\eqref{eq:powerflow_UC0} to ensure the positiveness of the variables:
\begin{subequations}
\label{eq:powerflow_UC}
\begin{align}
    & |p_{ik}(\tau)|\leq 2\bar{P}_{ik}, \forall \{i,k\}\in\E, \label{eq:powerflow_UC-1} \\
    & P(\tau) = \A (P_L(\tau) - \bar{P}_L), \label{eq:powerflow_UC-2} \\
    & P_L(\tau) = \B\A^\top\Theta(\tau) + \bar{P}_L, \label{eq:powerflow_UC-3}
\end{align}
\end{subequations}
where $\bar{P}_L$ is the collection of $\bar{P}_{ik}$ for all $\{i,k\}\in\E$.

All the constraints for the unit commitment problems have been included in~\eqref{eq:gen_bdds},~\eqref{eq:ramp_const},~\eqref{eq:OnOff_ind} and~\eqref{eq:powerflow_UC}. Our goal is to minimize the generation cost given as $\sum_{\tau=1,\cdots,T}\sum_{i\in\N_G}c_i p_i(\tau)$, where $c_i\in\real_+$ for all $i\in\N_G$. To this end, the complete formulation of the unit commitment problem is given as
\begin{align}\label{eq:unit_commit}
    &\min_{y\in\{0,1\}^{n_y},z\in\real^{n_z}_+} \sum_{\tau=1,\cdots,T}\sum_{i\in\N_G}c_i p_i(\tau) \\ \nonumber
    &\text{s.t. \eqref{eq:gen_bdds}, \eqref{eq:ramp_const}, \eqref{eq:OnOff_ind} and \eqref{eq:powerflow_UC},} 
\end{align}
where $y$ and $z$ are the collections of all the binary and continuous variables, respectively. Unit commitment problem~\eqref{eq:unit_commit} fits directly into the formulation of~\eqref{eq:p1g}. In fact, one can also choose to reformulate~\eqref{eq:unit_commit} to~\eqref{eq:p1} by compromising on the size of the problem. Specifically, we can introduce slack continuous variables to transform all the inequality constraints to equality constraints. The objective function of the continuous variables is approximated by binary variables with the method similar to~\eqref{eq:approx_q}. In this moderately sized test case, the additional slack variables do not pose major a challenge, so we choose to reformulate~\eqref{eq:unit_commit} to~\eqref{eq:p1} for both simpler coding and demonstration purposes. 
%(or \eqref{eq:p1} by introducing some slack variables). 
\begin{table}[H]
    \centering
    \begin{tabular}{|c|c|c|c|}\hline
        & G1 (bus 1) & G2 (bus 2) & G3 (bus 6) \\ \hline
        $c_1$ ($\$/$MW) & 16.83 & 40.62 & 21.93 \\ \hline
        %$c_0$ ($\$$) & 220.58 & 161.86 & 171.23 \\ \hline
        $\bar{P}$ (MW) & 200  & 100 & 20  \\ \hline
        $\underline{P}$ (MW) & 100 & 20 & 10 \\ \hline
        Ramp (MW/hr) & 50 & 40 & 15 \\ \hline
        $T_{on}$ (hr) & 2 & 2 & 2 \\ \hline
        $T_{off}$ (hr) & 2 & 2 & 2 \\ \hline
    \end{tabular}
    \caption{Parameters of the generators}
    \label{tab:6bus_gen}
\end{table}

\begin{table}[H]
    \centering
    \begin{tabular}{|c|c|c|c|c|} \hline
        Line & From bus & To bus & $X$ (p.u.) &  Flow limit (MW) \\ \hline
        1 & 1 & 2 & 0.17  & 200 \\ \hline
        2 & 1 & 4 & 0.258 & 100 \\ \hline
        3 & 2 & 4 & 0.197 & 100 \\ \hline
        4 & 5 & 6 & 0.14  & 100 \\ \hline
        5 & 3 & 6 & 0.018 & 100 \\ \hline
        6 & 2 & 3 & 0.037 & 100 \\ \hline
        7 & 4 & 5 & 0.037 & 100 \\ \hline
    \end{tabular}
    \caption{Parameters of the transmission lines}
    \label{tab:6bus_line}
\end{table}

\begin{table}[H]
    \centering
    \begin{tabular}{|c|c|c|c|}\hline
        Hour & L1 (bus 3) & L2 (bus 4) & L3 (bus 5) \\ \hline
        1 & 39.45 & 78.91 & 78.91 \\ \hline
        2 & 42.36 & 84.73 & 84.73 \\ \hline
        3 & 42.24 & 84.48 & 84.48 \\ \hline
    \end{tabular}
    \caption{Parameters of the loads (in MW)}
    \label{tab:6bus_load}
\end{table}

\rev{
\begin{table}[H]
    \centering
    \begin{tabular}{|c|c|c|c|c|}\hline
        $T$ & $n_y$ & $n_z$ & $m_b$ & $m_d$ \\ \hline
        2 & 87 & 68 & 30 & 72\\ \hline
        3 & 135 & 102 & 48 & 108 \\ \hline
    \end{tabular}
    \caption{Problem sizes of the optimization.}
    \label{tab:N_variables}
\end{table}
}

\subsubsection{Results with Binary Programming Solver}
\begin{figure}
    \centering
    \includegraphics[width=.48\textwidth]{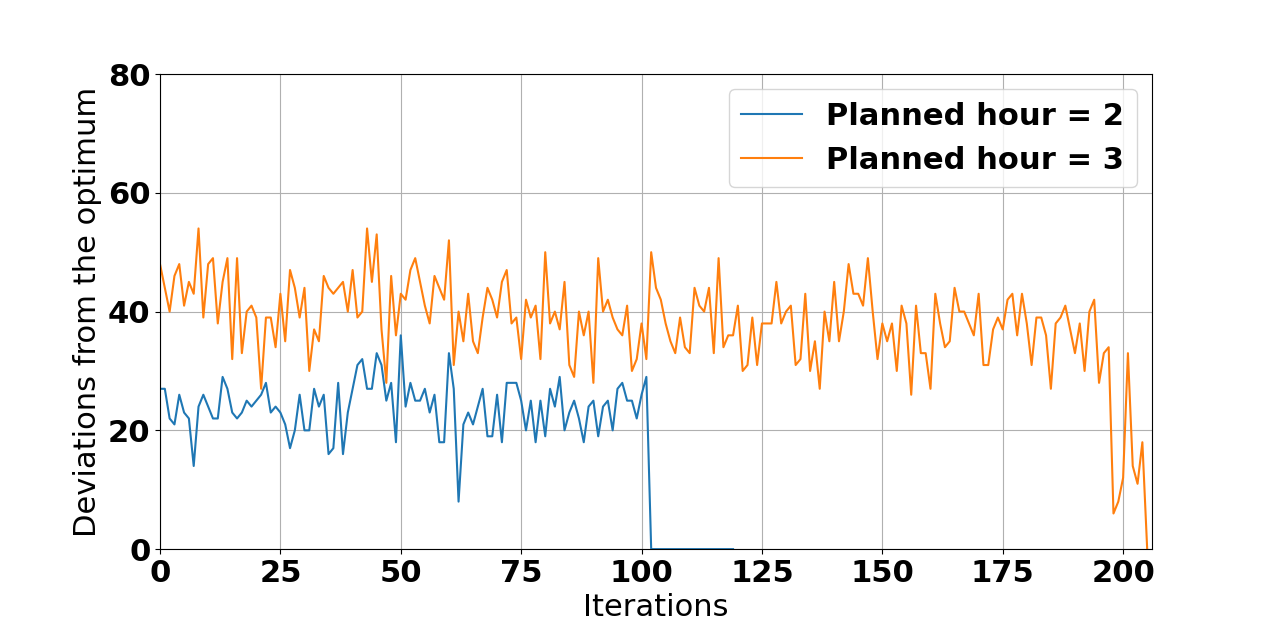}
    \caption{Convergence of Algorithm~\ref{algo:decomposition}. }
    \label{fig:Diff_to_opt_planned_3}
\end{figure}

\begin{figure}
    \centering
    \includegraphics[width=.48\textwidth]{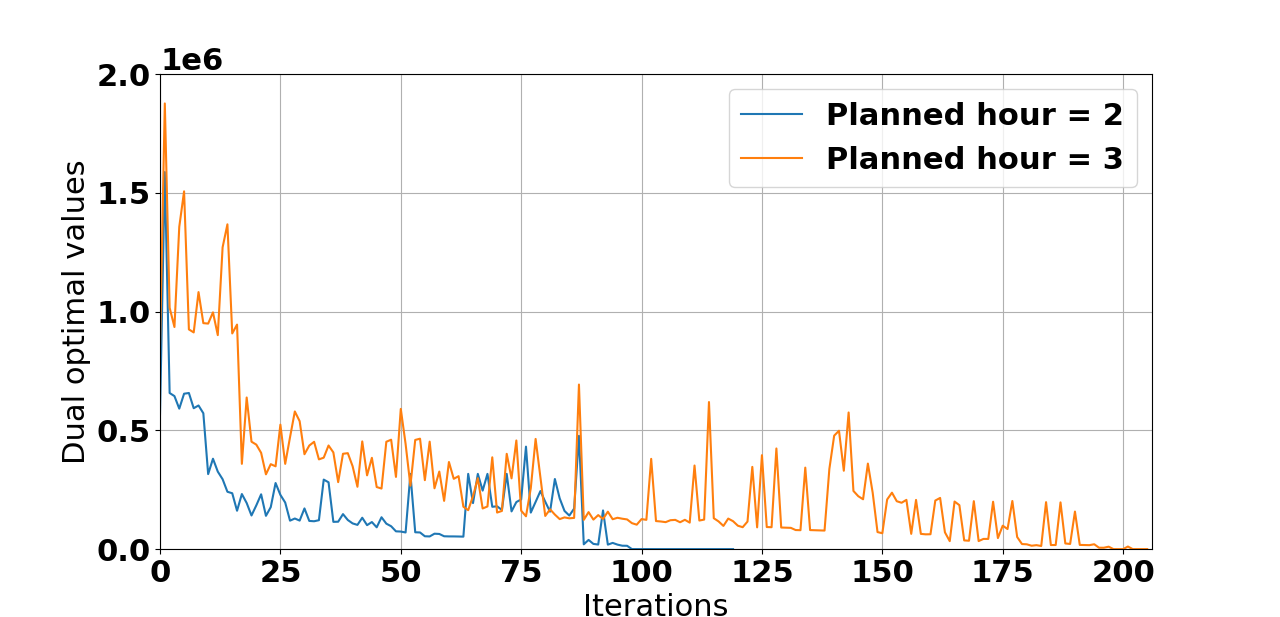}
    \caption{Evolution of the optimal values of~\eqref{eq:qfunction-d3}. }
    \label{fig:Dual_vals_planned_3}
\end{figure}

We first test Algorithm~\ref{algo:decomposition} against the unit commitment problem formulated in the previous subsection with the numerical setup detailed in Tables\rev{~\ref{tab:6bus_gen}-\ref{tab:N_variables}}. Both the binary programming and convex optimization of Algorithm~\ref{algo:decomposition} are solved by PuLP toolkit for python. The upper bound of every entry of $\alpha$, $\bar{\alpha}$, is set at 1000. We compare two scenarios with different planning time horizons (2 and 3 hours). As shown in Figure~\ref{fig:Diff_to_opt_planned_3}, Algorithm~\ref{algo:decomposition} finds the optimal solutions in both scenarios. The number of iterations it took to find the optimal solution grows with the problem size. The scenario with a 3 hour planning horizon has roughly 1.5 times more of the constraints and decision variables compared to the one with 2 hour, which results in roughly 100 more iterations (cuts) needed to reach the optimal solution. The simulation results give a rough picture on how the number of cuts grow with respect to the problem size. \revII{Figure~\ref{fig:Dual_vals_planned_3} shows that the optimal value of~\eqref{eq:qfunction-d3} decays gradually to zero as opposed to the abrupt jump to zero in Figure~\ref{fig:Diff_to_opt_planned_3}.} The optimal value of~\eqref{eq:qfunction-d3} serves as a good indicator of how close Algorithm~\ref{algo:decomposition} is to the optimum. 

\subsubsection{Simulated Annealing Sampler} Simulated annealing algorithm is based on the technique of cooling metal from a high temperature to improve its structure (annealing)~\cite{Dwave_neal}. We choose a simulated annealing sampler for the unit commitment problem largely due to hardware constraints of current D-Wave systems. Results in the last subsection show that Algorithm~\ref{algo:decomposition} needs to add \rev{104} constraints (cuts) before reaching the optimal solution. Each additional inequality constraint requires introducing concomitant slack binary variables as discussed in Section~\ref{sec:reformulation}. The number of variables added depends on the desired floating point precision. Even if we choose a relatively large floating point precision at 1 (granularity constant $w_{x_k}=1$ in~\eqref{eq:binary_qc1}), around 15 additional binary variables are needed for each new constraint. The binary program in Algorithm~\ref{algo:decomposition} will therefore grow to a scale of around 2000 binary variables in the last few iterations. Such binary programs are not solvable by the D-Wave 2000Q system that we have access to because embedding the binary program in the D-Wave hardware graph requires physical qubit overhead that results in far more than the 2000 logical qubits needed for the 2000 binary variable problem. The simulated annealing sampler follows similar heuristics to the D-Wave system although it is run on classical computers and relies on ``thermal'' hopping rather than quantum tunneling to traverse a problem's energy landscape. We therefore view simulated annealing as a good proxy for an actual D-Wave processor for algorithm validation purposes. \rev{The view aligns with the literature well~\cite{karimi2017effective,irie2021hybrid}.}  

We compare the performance of Algorithm~\ref{algo:decomposition} with two different implementations: one uses PuLP to solve the binary programs involved in Algorithm~\ref{algo:decomposition} as shown in the last subsection; the other uses simulated the annealing sampler. Both implementations are tested on the unit commitment problem with $T = 2$. For simulated annealing, we set the penalty strength at 2 (one need not set the chain strength because it is only required for an actual D-Wave system). 
%Different from actual D-Wave systems, simulating annealing does not need to tune the chain strength. 
In addition to the penalty strength, other parameters, including number of reads ($N_{read}$) and number of sweeps ($N_{sweep}$), have a significant impact on simulated annealing results. $N_{read}$ is the number of runs of the simulated annealing algorithm (\rev{analogous to the number of runs of quantum annealing}). $N_{sweep}$ is the number of fixed temperature sampling sweeps \rev{for each run of the simulated annealing algorithm}. As the size of the optimization problem grows, one requires larger values of those parameters for similar precision of the solutions. Therefore, we take the approach of gradually increasing $N_{read}$ and $N_{sweep}$ as the number of iterations increases in Algorithm~\ref{algo:decomposition}. Specifically, we set $N_{read} = 80$ and $N_{sweep}= 5000$ at $t=0$ and allow them to grow exponentially with the rate of $0.02$, given as
\begin{align}\label{eq:exp_grow}
    &N_{read}(t) = \lfloor 80\cdot(1.02)^t\rfloor,\\ \nonumber  
    &N_{sweep}(t) = \lfloor 5000\cdot (1.02)^t \rfloor,
\end{align}
where $t$ is the iteration number of the algorithm. 

\begin{figure}
    \centering
    \includegraphics[width=.48\textwidth]{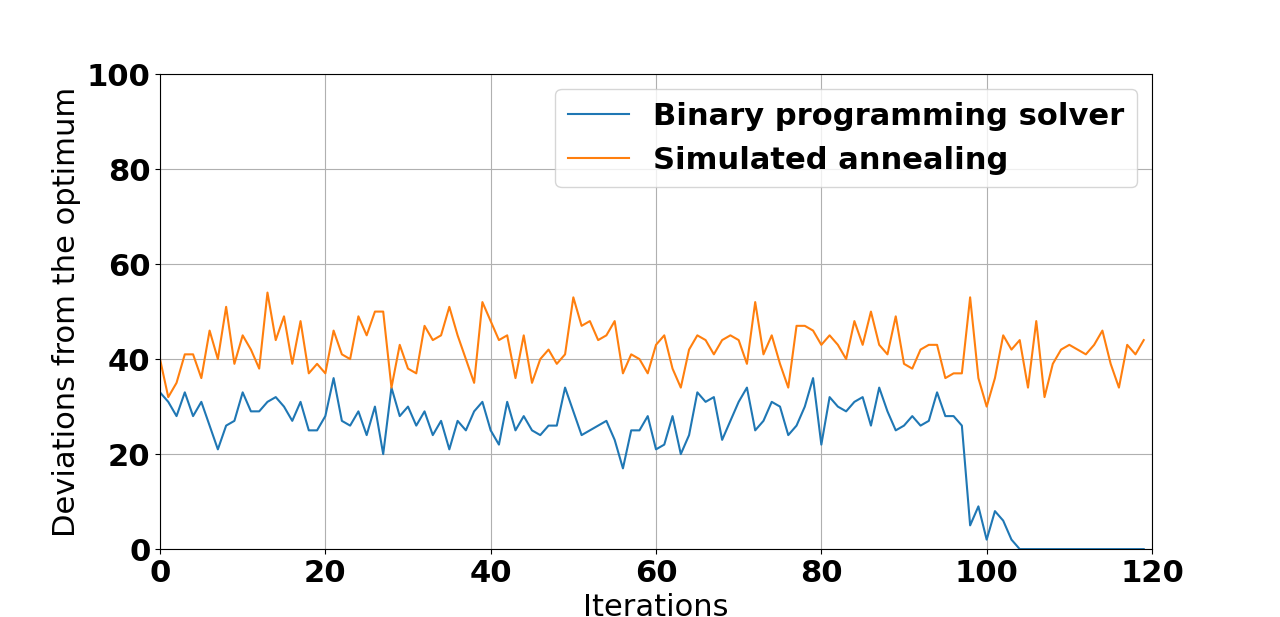}
    \caption{Convergence to the optimal point}
    \label{fig:Diff_to_opt}
\end{figure}

\begin{figure}
    \centering
    \includegraphics[width=.48\textwidth]{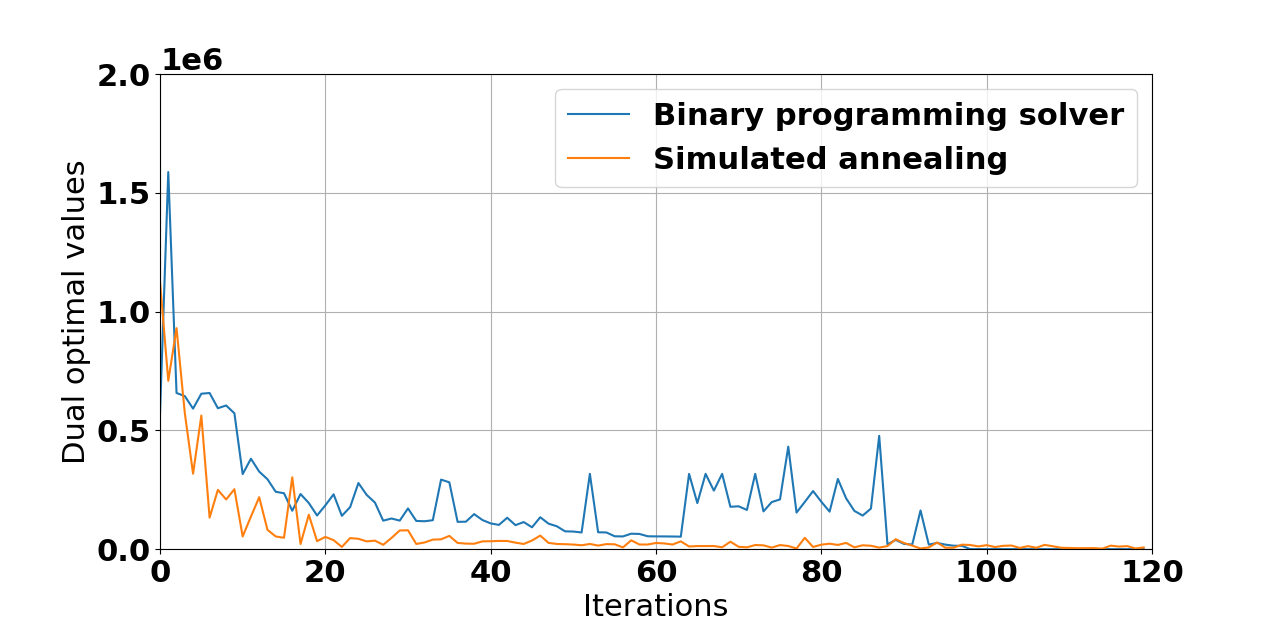}
    \caption{Evolution of the optimal values of~\eqref{eq:qfunction-d3}. }
    \label{fig:Dual_vals}
\end{figure}

\begin{figure}
    \centering
    \includegraphics[width=.48\textwidth]{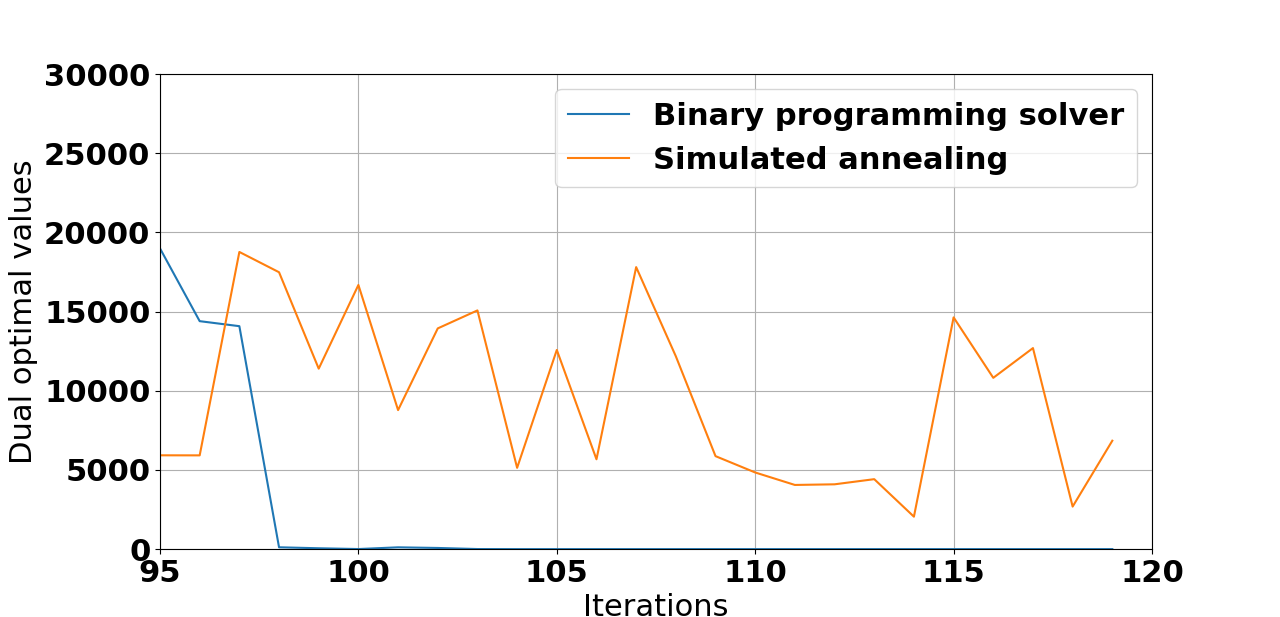}
    \caption{Evolution of the optimal values of~\eqref{eq:qfunction-d3} (zoomed in).}
    \label{fig:Dual_vals_2}
\end{figure}

\begin{figure}
    \centering
    \includegraphics[width=.48\textwidth]{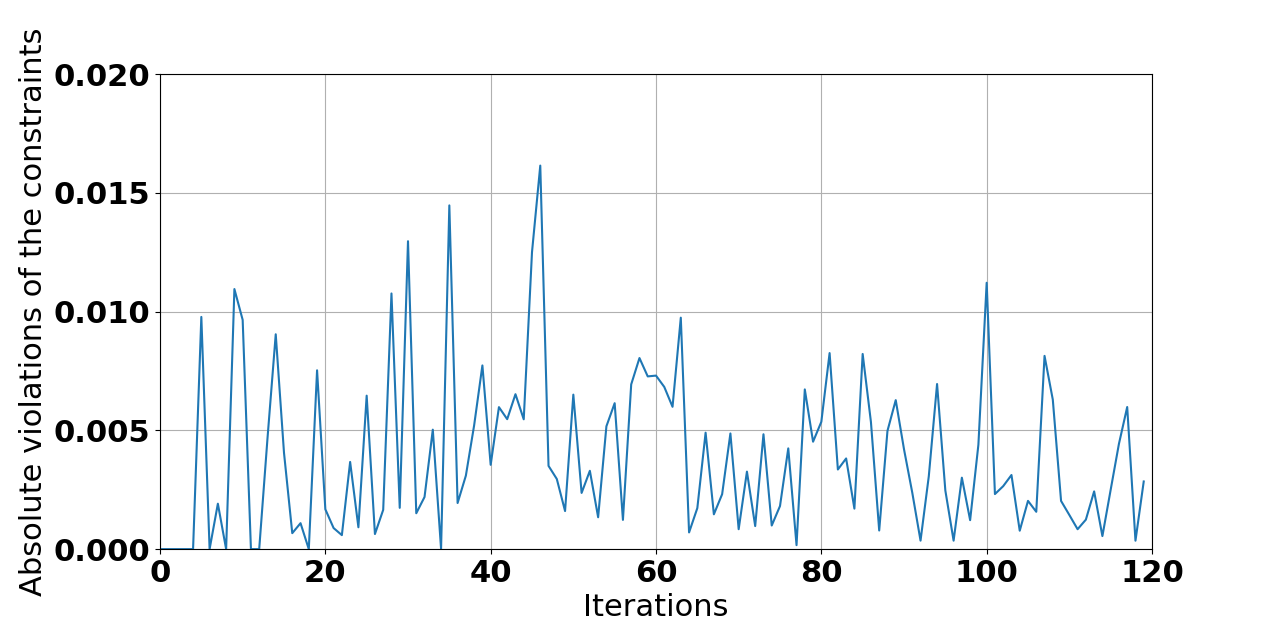}
    \caption{Violations of the constraints for simulated annealing sampler. }
    \label{fig:tot_deviations_normalized}
\end{figure}

As shown in Figure~\ref{fig:Diff_to_opt}, PuLP obtains the optimal solution in 104 steps while simulated annealing fails to find the optimum with 120 steps. \rev{The optimal solution that PuLP found has G1 and G2 in the ``on'' state and G3 in the ``off'' state for all time. The solution that simulated annealing reached has all the generators in the ``on'' states. Both solutions are feasible with a marginal difference between the optimal values, which are respectively 5986 and 6021. Simulated annealing version still found a reasonably good solution with less than $1\%$ off the optimum. A similar conclusion can be drawn by observing the evolution of the optimal values of~\eqref{eq:qfunction-d3} in Figures~\ref{fig:Dual_vals} and~\ref{fig:Dual_vals_2}, which show that the simulated annealing approaches the optimal solution as the number of iterations increase.
}
The main reason that simulated annealing fails to find the optimal solution is the violation of the constraints (cuts). Recall that in simulated annealing, we penalize the inequality constraints by quadratic functions given by~\eqref{eq:binary_qc1-2}. We quantify the violations by the following
\begin{align}\label{eq:vio_abs}
    \sum_{k = 1}^{t} \frac{1}{t\cdot \hat{b}_t }\Big|a_{k} y - \hat{b}_t(k) + w_{x_k} \sum_{j = 1}^{n_{x_k}}2^{(j-1)}x_k(j)\Big|.
\end{align}
%Figure~\ref{fig:tot_deviations_normalized} shows the value of~\eqref{eq:vio_abs} as a function of iteration, $t$. 
\rev{\revII{The violations as a function of $t$ shown in Figure~\ref{fig:tot_deviations_normalized} are within a small margin but fail to drop further to the exact optimal solution as the PuLP case.} This reflects one of the challenges as the number of constraints (cuts) grows. Adding a constraint in the form of an extra quadratic penalty function effectively diminishes the weight share of the previous penalty functions and the original objective function. This is not an issue from the analytical standpoint because as long as the penalty weights are non-zero, the QUBO always has the same optima as the original constrained binary programming as discussed in Section~\ref{sec:reformulation}. However, numerically, especially for simulated annealing or QPU that rely heavily on the probabilistic distribution of samples, it is challenging to find penalty weights that guarantee distinction between the global optima from near optimal solutions. We have many trials of various penalty weights, but they all generate similar results in the sense that reaching near global optima in the last few iterations. Using the metric in equation~\eqref{eq:qfunction-d3} (illustrated in Figure~\ref{fig:Dual_vals_2}) for early termination on a near optimum is a reasonable way to implement Algorithms~\ref{algo:decomposition} and~\ref{algo:decomposition-2}. 
}

%Despite the fact that the violations are within a rather small margin, simulation results show that for Algorithm~\ref{algo:decomposition} to converge on quantum systems, one needs to carefully tune the penalty strengths such that the violations are small enough to the point where the added cuts can effectively steer the solutions to the optimal one.

\rev{We last want to comment on the runtime for the simulated annealing. Each emulation of a result from quantum annealer needs $N_{read}\cdot N_{sweep}$ number of samples, which translates to the order of $10^7$ (or more) samples from iteration 100 onward. Therefore, the runtime for generating Figure~\ref{fig:tot_deviations_normalized} is around 37 hours (did not time exactly) on a personal laptop with 3.5GHz CPU and 8GB RAM. This is a lot longer than many existing solvers. However, simulated annealing sampler only serves the purpose of generating similar results of quantum annealer. We expect the actual runtime on quantum annealer to be a lot faster, and advantages over MIP standard solvers could be observed in large-scale MIP in future endeavors.}

%The fact that the average violation does not grow with respect to the $t$ reflect that~\eqref{eq:exp_grow} has $N_{read}(\cdot)$ and $N_{sweep}(\cdot)$ grow at a rate that matches the increasing size of the optimization. 

\section{Conclusion}\label{sec:conclusion}
In this paper, we proposed applying Benders decomposition for a hybrid QC-CC algorithm. The algorithms presented take an iterative approach that gradually closes the gap the optimal solution. Small-scale simulation results demonstrate that the algorithms can find the global optimum with proper selection of D-Wave quantum annealer parameters. For a more realistic MIP that is inspired by unit commitment applications, we identify some limitations of the effectiveness of the algorithm mainly due to the random nature of the annealing mechanism and the requirement that hard constraints be softened in such approaches. \rev{In addition to the challenges of the random nature and methods of softening the constraints, another important improvement need of the proposed algorithm is characterizing the bound of the number of iterations required for finding a reasonably good solution. 
Our near future work will focus more on some tangible tasks, including rigorous methods for identifying penalty weights for robust numerical performance of the hybrid algorithm. We will also look towards additional real-world MIP applications and refine the algorithms for those applications.}

%\section*{Acknowledgement}
%We like to thank Amazon Web Services (AWS) for providing access to the 2000Q D-Wave system for our initial testing on the algorithm before their official launch. 

	\bibliographystyle{IEEEtran}
	\bibliography{ref}
\end{document}